# Chi-square Intervals for a Poisson Parameter
# — Bayes, Classical and Structural —

E.A. Maxwell*

**Abstract.** The 'standard' confidence interval for a Poisson parameter is only one of a number of estimation intervals based on the chi-square distribution that may be used in the estimation of the mean or mean rate for a Poisson model. Other chi-square intervals are available for experimenters using Bayesian or structural inference methods. Exploring these intervals also leads to other alternate approximate chi-square intervals. Although coverage probability may not always be of interest for Bayesian or structural intervals, coverage probabilities are useful for validating 'objective' priors. Coverage probabilities are explored for all of the intervals considered.

**MSC-class:** 62F25; 62F15, 62A01

## 1 Introduction

Recent interest in interval estimation for the mean of a Poisson distribution has included methods to revise/improve the standard chi-square interval to produce an interval that is optimal in the sense of being shorter than other possible intervals while maintaining a coverage probability that is at least equal to the nominal confidence level. The methods require that the experimenter run a computer algorithm to obtain the resulting interval. (Casella and Robert, 1989; Wardell, 1997; Kabaila and Byrne, 2001.)

These methods are intended to produce intervals that are exact in the sense that any resulting $100(1 - \alpha)\%$ confidence interval will have a coverage probability that is always at least $1 - \alpha$ (Casella and Robert, 1989; Kabaila and Byrne, 2001,) although that goal is not always achieved (Berger and Coutant, 2001.) They are not "exact" in the way that the usual confidence interval for the mean of a normal population is exact — that is, they do not have coverage probability of $(1 - \alpha)$ in every case. It may be argued (Agresti and Coull, 1998) that it is better to use approximate intervals rather than "exact" ones and to obtain a coverage probability that is close to $(1-\alpha)$ rather than one that is always at least $(1-\alpha)$ but often much larger than $(1-\alpha)$.

*Department of Mathermatics, Trent University, Peterborough, ON, Canada  cheam2@sympatico.ca

Developing and running a computer algorithm, including a spreadsheet algorithm, may appeal to many analysts; on the other hand, many investigators may prefer simpler methods such as those based on basic distribution percentiles available from tables or internally in computer software/spreadsheets. A number of chi-square-based intervals are explored below. These include the usual interval which may be longer than necessary and have a coverage probability greater than (1-$\alpha$) plus several shorter intervals with coverage probabilities that can be less than (1-$\alpha$), but may generally be close to (1-$\alpha$). The shorter intervals studied below include intervals that are derived not as confidence intervals, but as structural probability intervals or Bayes credible intervals. Two additional intervals suggested by these are explored as well.

Because structural inference and Bayesian inference produce post-sample probability distributions for the parameter(s) of interest, the resulting estimation intervals are interpreted differently from confidence intervals and have exact posterior probability (1-$\alpha$) of capturing the parameter of interest. Coverage probabilities are not necessarily relevant or of interest to the investigator using such intervals. It has been suggested, though, (e.g. Ghosh and Ghosh, 2004) that coverage probabilities for Bayes credible intervals may be used for validating nonsubjective priors. Coverage probabilities are explored for all of the intervals considered.

The following notation is used. A chi-square variable with $f$ degrees of freedom is denoted as $\chi^2_{[f]}$. Its $100p^{\text{th}}$ percentile is denoted as $\chi^2_{[f],p}$ ($0 \leq p \leq 1$) and is such that

$$P\left[\chi^2_{[f]} \leq \chi^2_{[f],p}\right] = p$$

In each case, $\nu$ is the mean rate for the Poisson process under investigation, $\lambda$ is the mean number of occurrences for the time/distance/area/volume sampled and $t$ is the time/distance/area/volume sampled so that $\lambda = \nu t$. The observed number of occurrences in the sample is denoted as $x$.

The resulting 100(1 - $\alpha$)% estimation intervals for $\lambda$ all have the form:

$$\frac{1}{2} \chi^2_{[f_1]\alpha/2} \quad to \quad \frac{1}{2} \chi^2_{[f_2](1-\alpha/2)}$$

Where $f_1$ and $f_2$ are functions of $x$.

The resulting 100(1 - $\alpha$)% estimation intervals for $\nu$ then all have the form:

$$\frac{1}{2t} \chi^2_{[f_1]\alpha/2} \quad to \quad \frac{1}{2t} \chi^2_{[f_2](1-\alpha/2)}$$

The case of $x = 0$ produces difficulties. To determine an interval for each case for comparison purposes, the following procedure was used. In the event that $f_1 = 0$ (e.g. $f_1 = 2x$ and $x = 0$,) the lower confidence limit is taken to be 0. In the event that $f_2 = 0$ (e.g. $f_2 = 2x$ and $x = 0$,) the upper confidence limit is calculated with $f_2 = 1$.

For each of the cases of $f_1$ and $f_2$ that arise from the different approaches explored, both 95% and 99% intervals are considered. Coverage probabilities were determined for $\lambda = 0.1$ to 75.0 in steps of 0.1.

## 2 "Usual" intervals

The "usual" chi-square interval based on the relationship between the Poisson and chi-square distribution functions (e.g. Johnson & Kotz, 1969 Chap 4) has $f_1 = 2x$ and $f_2 = 2x+2$. Its coverage probability has already been well studied and plotted; however, for completeness it is included again below. As illustrated in Figures 1 and 2, the coverage probability is usually, and sometimes substantially, larger than $(1-\alpha)$, indicating that the intervals tend to be longer than necessary, thus, generating interest in alternative, shorter intervals.

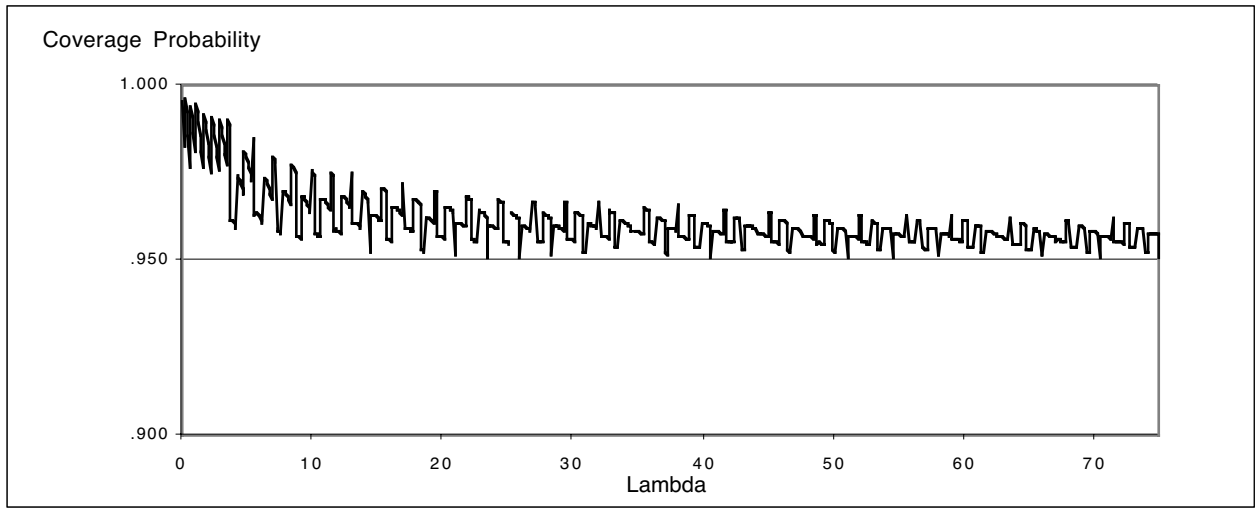

FIGURE 1: Coverage Probabilities for 95% Intervals for $f_1 = 2x$ and $f_2 = 2x + 2$

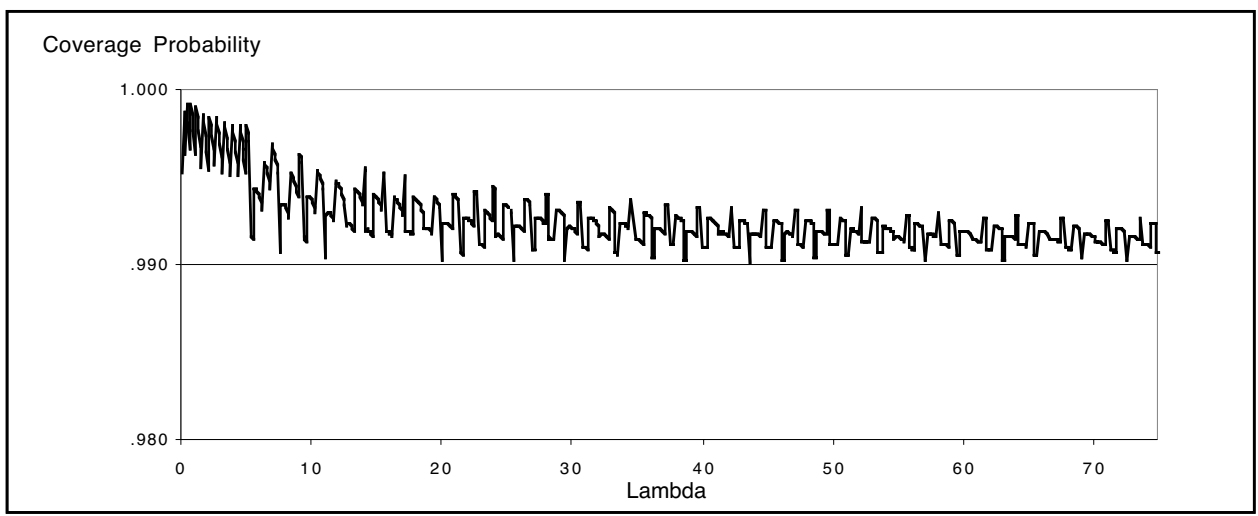

FIGURE 2: Coverage Probabilities for 99% Intervals for $f_1 = 2x$ and $f_2 = 2x + 2$

## 3  Structural/Bayes intervals

Because it involves a discrete random variable, the Poisson model does not lend itself directly to "standard" procedures of structural inference (Fraser, 1968.)  Structural inference procedures can be applied indirectly, though, in a combination structural/Bayes approach.  Interarrival times, distances, etc. or waiting times, distances etc. until the $n^{\text{th}}$ Poisson event do fit in the framework of the multiplicative error model or scale model (Fraser, 1968, Chap. 1.)

For example, for events occurring in time at a mean rate $\nu$, interarrival times follow the exponential distribution with mean $\theta = 1/\nu$.  If $n$ sample interarrival times $t_1$, $t_2$, $\cdots t_n$ produce a total time $t$, then the density function resulting from the structural distribution for $\theta$ based on $t$ is proportional to

$$\frac{t^n}{\theta^{n+1}} e^{-t/\theta}$$

and the resulting density function for $\nu$ based on $t$ is

$$f\left(\nu \mid t\right) = \frac{t^n}{\Gamma\left(n\right)} \nu^{n-1} e^{-\nu t}$$

i.e. the structural distribution for $\nu$ based on $t$ is a gamma distribution with parameters $n$ and $1/t$. The result is the same if a process is observed until the $n^{\text{th}}$ occurrence and $t$ is the waiting time until the $n^{\text{th}}$ occurrence.

This result may be used as part of a tandem design similar to that suggested for using structural inference in investigations involving type I censoring (Whitney and Minder 1974.) The tandem design produces a distribution for the parameter of interest which Whitney and Minder note is not a usual structural distribution, but actually a Bayes posterior distribution developed from a structural distribution as a prior from the first stage.

In the first stage, a process is observed until the $n_1^{\text{th}}$ occurrence and the waiting time $t_1$ noted. The resulting structural density for $\nu$ is

$$f\left(\nu \mid t_1\right) = \frac{t_1^{n_1}}{\Gamma\left(n_1\right)} \nu^{n_1-1} e^{-\nu t_1}$$

This is used as a prior density for the second stage.  In the second stage, the process is observed for a further time $t_2$ and the number of occurrences $n_2$ is noted.  The likelihood function is proportional to

$$\left(\nu t_2\right)^{n_2} e^{-\nu t_2} / n_2!$$

and the resulting posterior for $\nu$ then is proportional to

$$\nu^{(n_1+n_2-1)} e^{-\nu(t_1+t_2)}$$

This is a function of the total number of occurrences $n_1 + n_2$ which may be written as $x$ and the total time observed $t_1 + t_2$ which may be written as $t$. After normalizing, the structural-inference-based posterior density for $\nu$ based on $t$ then is

$$f(\nu \mid t) = \frac{t^x}{\Gamma(x)} \nu^{x-1} e^{-\nu t}$$

i.e. the structural-inference-based posterior for $\nu$ is a gamma distribution with parameters $x$ and $1/t$. This distribution is independent of the individual values for $n_1$ and $n_2$; it depends only on their total $x$. Similarly it is independent of the individual values for $t_1$ and $t_2$, depending only on their total $t$. To avoid an excessively large value for $t_1$, $n_1$ should be small. Indeed, it is possible to choose $n_1 = 1$.

With $\lambda = \nu\, t$, the structural-inference-based posterior for $\lambda$ is a gamma distribution with parameters $x$ and 1. The structural-inference-based posterior for $2\lambda$ then is a gamma distribution with parameters $x$ and 2, i.e. it is a chi-square distribution with $2x$ degrees of freedom. This result provides a 100(1-$\alpha$)% interval for $2\lambda$ and, hence, a 100(1-$\alpha$)% chi-square interval for $\lambda$ with $f_1 = 2x$ and $f_2 = 2x$.

$$\frac{1}{2}\chi^2_{[2x]\alpha/2} \quad to \quad \frac{1}{2}\chi^2_{[2x](1-\alpha/2)}$$

The resulting 100(1 - $\alpha$)% interval for $\nu$ is

$$\frac{1}{2t}\chi^2_{[2x]\alpha/2} \quad to \quad \frac{1}{2t}\chi^2_{[2x](1-\alpha/2)}$$

It is not necessary to restrict this form of analysis to a model of events in time. The above method will work with any linear medium; e.g. $t$ can be distance rather than time. For events in the plane or in space, the sampling method can still be used. An initial point is selected in the plane or space and then the surrounding area or volume is expanded radially (within the boundaries of the available sample space) until the $n_1$th event is observed and the area or volume sampled, $t_1$, is noted. A further area or volume $t_2$ is then sampled and the number of events observed, $n_2$, is noted. The results then are as above.

With $f_2 = 2x$ instead of $2x+2$, these intervals are shorter than the usual intervals and have lower coverage probabilities. In calculating these coverage probabilities, it is interesting to note that, when taking into account the probability models related to the two stages of the tandem design, and using the combined probability models, the ultimate probability of obtaining a total of $x$ occurrences in a total time $t$ is

$$\lambda^x e^{-\lambda} / x!$$

where $\lambda = \nu t$; i.e. it is just the basic Poisson probability of obtaining a total of $x$ occurrences in a total time. $t$.

Coverage probabilities are illustrated in Figures 3 and 4.

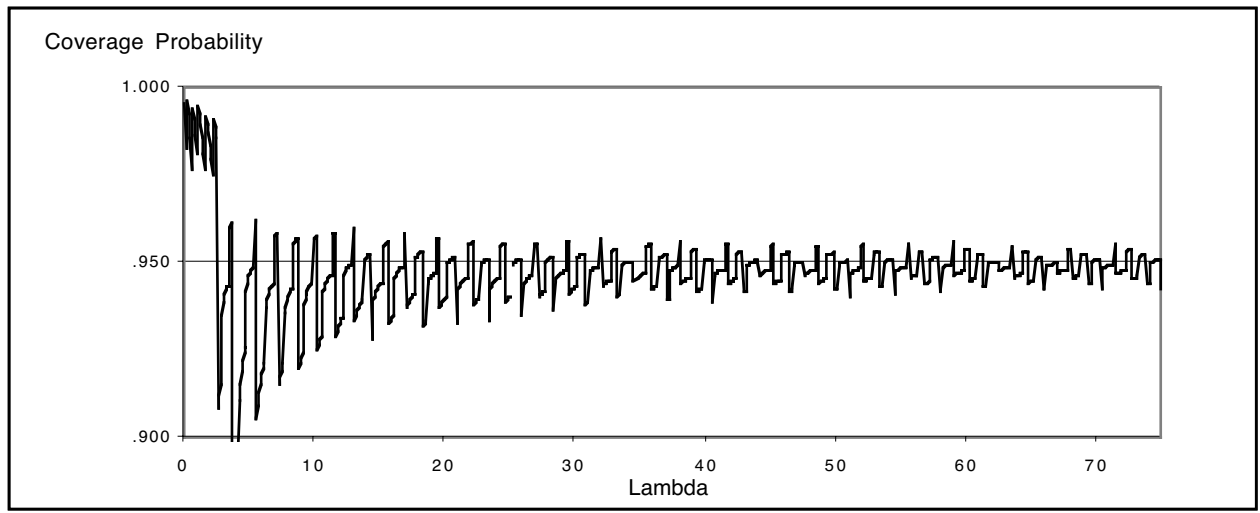

FIGURE 3: Coverage Probabilities for 95% Intervals for $f_1 = 2x$ and $f_2 = 2x$

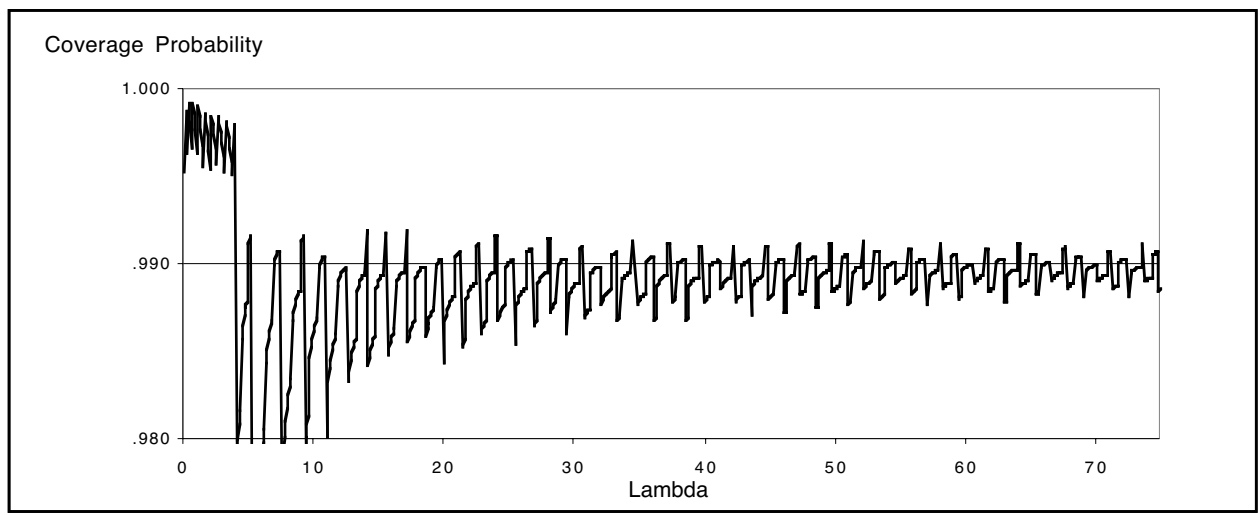

FIGURE 4: Coverage Probabilities for 99% Intervals for $f_1 = 2x$ and $f_2 = 2x$

Figures 3 and 4 indicate that, for these intervals, the coverage probabilities tend to be much smaller than the nominal ($1-\alpha$) level. For the analyst comfortable with this structural approach, that is not a problem. The intervals have a different interpretation; they have a structural-probability-based post-sample probability ($1-\alpha$) of including the parameter of interest.

For investigators for whom the coverage probability is important, these intervals would appear to be too short and the coverage probabilities too low.

# 4. Bayesian intervals

Bayesian inference models are explored within the framework of conjugate priors (e.g. Raiffa and Schlaifer, 1968; Lee, 1989; Press, 1989.) It is assumed that the investigator is comfortable working with improper priors. With $\lambda = v\, t$, the likelihood for a total of $\Sigma x$ Poisson events in a sample of $n$ times/distances/areas/volumes $t$ is proportional to

$$\lambda^{\Sigma x} e^{-n\lambda} / (\Pi x)!$$

and so a conjugate prior density for $\lambda$ is proportional to

$$\lambda^{a} e^{-b\lambda}$$

Using such a prior and letting $x$ represent the total number of observed Poisson events in one sample with a total time/distance/area/volume $t$, the posterior density for $\lambda$ is proportional to

$$\lambda^{a+x} e^{-(b+1)\lambda}$$

With appropriate hyperparameters $a$, $b$, the prior and posterior distributions are gamma distributions and the posterior distribution for $2\lambda$ becomes a chi-square distribution. The chi-square distribution then produces posterior estimation intervals or Bayesian credible intervals or credibility intervals (e.g. Press, 1989.)

Three priors which may be considered to be free of subjective influence are explored — a uniform prior, a Jeffreys prior for a scale parameter, and a Jeffreys prior particular to the Poisson model.

The uniform prior for $\lambda$ corresponds to $a = 0$ and $b = 0$ and produces a posterior for $\lambda$ that is proportional to

$$\lambda^{x} e^{-\lambda}$$

and this leads to the posterior density for $2\lambda$ as

$$f\left(2\lambda \mid x\right) = \frac{1}{2^{\left(\frac{2x+2}{2}\right)} \Gamma\left(\frac{2x+2}{2}\right)} (2\lambda)^{\left(\frac{2x+2}{2}-1\right)} e^{-2\lambda/2}$$

i.e. the posterior distribution for $2\lambda$ is chi-square with $2x + 2$ degrees of freedom. A Bayes credible interval for this case, thus, is again a chi-square interval, but with $f_1 = f_2 = 2x + 2$.

This interval is shorter than the usual chi-square confidence interval because $f_1$ has been increased from $2x$ to $2x + 2$. As a result it will have lower coverage probabilities.

Coverage probabilities are illustrated in Figures 5 and 6.

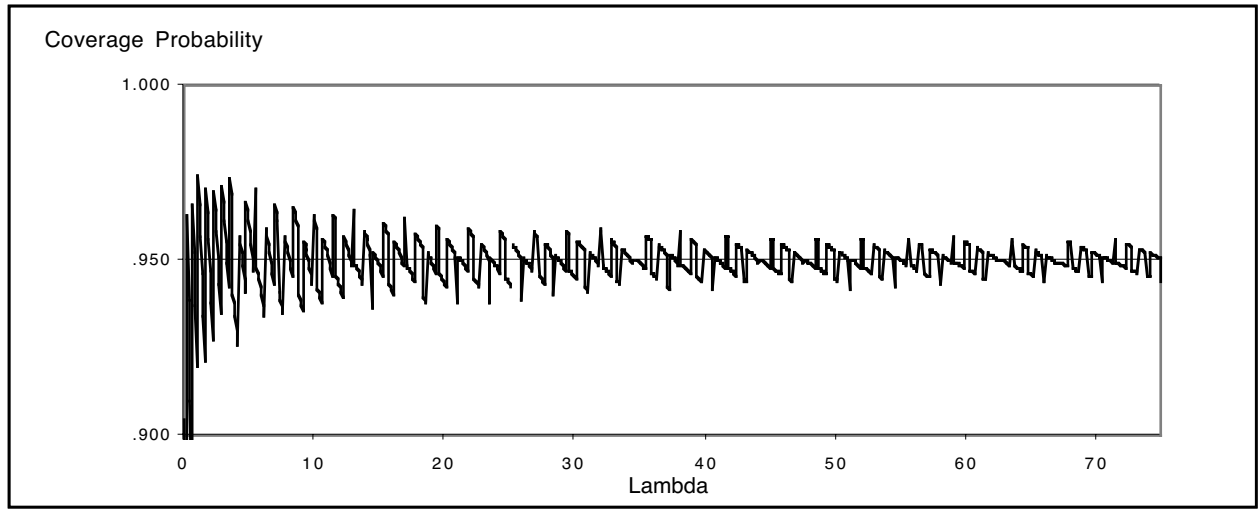

FIGURE 5: Coverage Probabilities for 95% Intervals for $f_1 = 2x + 2$ and $f_2 = 2x + 2$

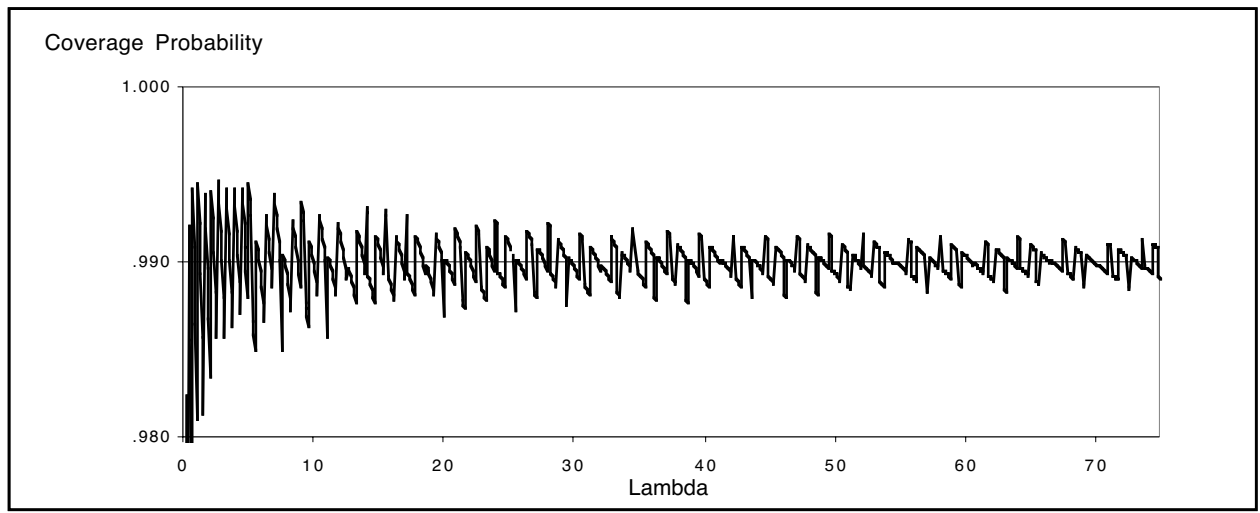

FIGURE 6: Coverage Probabilities for 99% Intervals for $f_1 = 2x + 2$ and $f_2 = 2x + 2$

Figures 5 and 6 indicate that the coverage probability can be quite low for cases with small $\lambda$, but the probabilities do fluctuate reasonably around $1-\alpha$ for other values of $\lambda$. As with the structural-inference-based intervals, the Bayes intervals have a different interpretation from that for confidence intervals. For the investigator comfortable with the Bayes approach and the chosen prior, the intervals will have post-sample probability of $1 - \alpha$ of including $\lambda$.

For investigators for whom the coverage probability is important, these intervals would appear to be reasonable if $\lambda$ is not too small, especially for investigators who prefer approximate intervals rather than "exact" ones to obtain a coverage probability that is close to $(1-\alpha)$ rather than one that is always at least $(1-\alpha)$ but often much larger than $(1-\alpha)$.

For a positive parameter $\lambda$, Jeffreys suggested a prior proportional to $1/\lambda$ (Jeffreys, 1961, Chap 3; Lee, 1989, Chap. 3.) This prior for $\lambda$ corresponds to $a = -1$ and $b = 0$ and produces a posterior distribution for $2\lambda$ which is chi-square with $2x$ degrees of freedom. A Bayes credible interval for this case, thus, is a chi-square interval with $f_1 = f_2 = 2x$. This is the same as the structural-inference-based interval with coverage probabilities illustrated in Figures 3 and 4.

For the particular case of the Poisson with parameter $\lambda$, the Jeffreys' prior is proportional to $1/\sqrt{\lambda}$ (Jeffreys, 1961, Chap 3; Lee, 1989, Chap. 3.) This prior for $\lambda$ corresponds to $a = -1/2$ and $b = 0$ and produces a posterior distribution for $2\lambda$ which is chi-square with $2x + 1$ degrees of freedom. A Bayes credible interval for this case is a chi-square interval with $f_1 = f_2 = 2x + 1$.

This interval is shorter than the usual chi-square confidence interval because $f_1$ has been increased from $2x$ to $2x + 1$ and $f_2$ has been decreased from $2x + 2$ to $2x + 1$. As a result this interval will again have lower coverage probabilities.

Coverage probabilities are illustrated in Figures 7 and 8. The coverage probability can be quite low for a very few cases with small $\lambda$; otherwise, the probabilities fluctuate quite well around $1 - \alpha$ for other values of $\lambda$.

Again, these are Bayes credible intervals and have a different interpretation from that for confidence intervals. For the investigator comfortable with the Bayes approach and the chosen prior, the intervals will have post-sample probability of $1 - \alpha$ of including $\lambda$.

For investigators for whom the coverage probability is important, these intervals would appear to be to reasonable except for a very few cases in which $\lambda$ is small. As above, this is especially so for investigators who prefer approximate intervals rather than "exact" ones to obtain a coverage probability that is close to $(1-\alpha)$ rather than one that is always at least $(1-\alpha)$ but often much larger than $(1-\alpha)$.

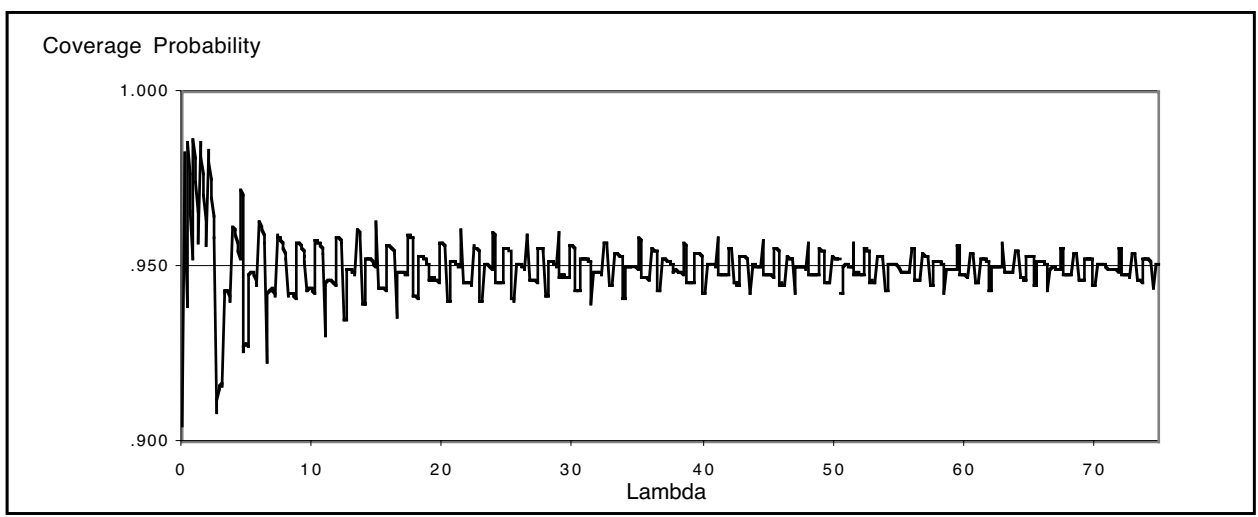

FIGURE 7: Coverage Probabilities for 95% Intervals for $f_1 = 2x + 1$ and $f_2 = 2x + 1$

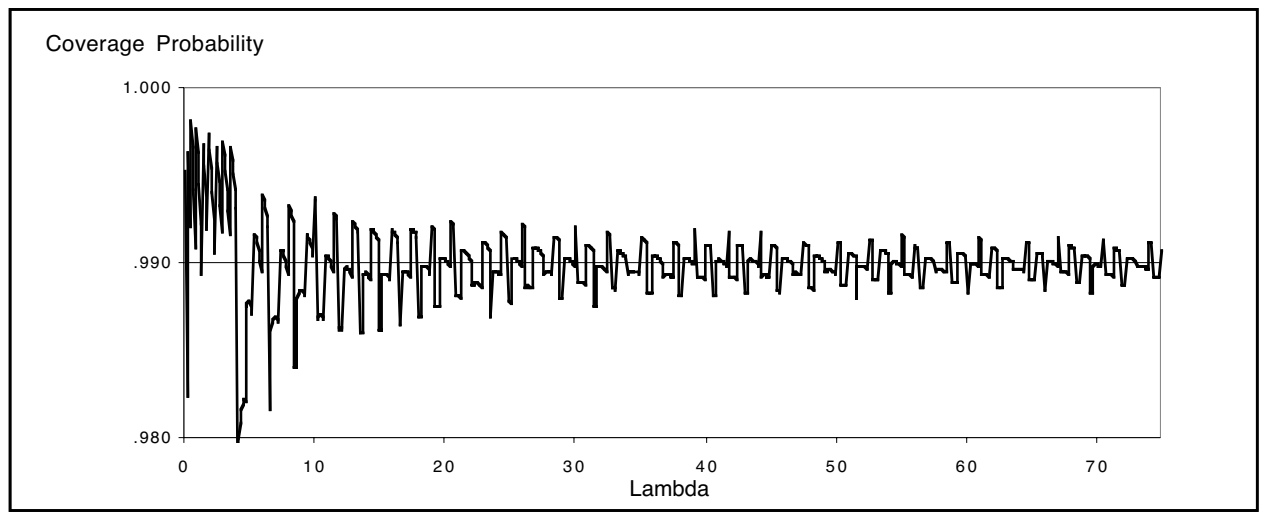

FIGURE 8: Coverage Probabilities for 99% Intervals for $f_1 = 2x + 1$ and $f_2 = 2x + 1$

That these intervals (with $f_1 = f_2 = 2x + 1$) have reasonable coverage probabilities may be taken by some investigators as providing validation of the Jeffreys' prior as a reasonably objective prior for the Poisson model.

## 5. Other possibilities

In each of the intervals above, the width of the interval is shorter than that for the "usual" interval due to what may be called a two-degrees-of-freedom adjustment. In the structural-inference-based interval or the Bayes interval based on a prior proportional to $1/\lambda$, $f_2$ is decreased by two degrees of freedom from $2x + 2$ to $2x$. For the Bayes interval based on a uniform prior, $f_1$ is increased by two degrees of freedom from $2x$ to $2x + 2$. For the Bayes interval based on the Jeffreys prior for the Poisson, $f_1$ is increased by one degree of freedom from $2x$ to $2x + 1$ and $f_2$ is decreased by one degree of freedom from $2x + 2$ to $2x + 1$.

The two-degrees-of-freedom adjustments produce shorter intervals with lower coverage probabilities — perhaps too low for some investigators. A possible compromise between these new intervals (which do have some theoretical justification and alternate interpretation) is a one-degree-of-freedom adjustment; i.e. increase $f_1$ by one degree of freedom from $2x$ to $2x + 1$ or decrease $f_2$ by one degree of freedom from $2x + 2$ to $2x + 1$. The motivation for such an adjustment is simply an exploratory trial of intervals similar to those developed above for the sake of checking coverage probabilities.

Coverage probabilities for these two other possibilities are illustrated in Figures 9 through 12.

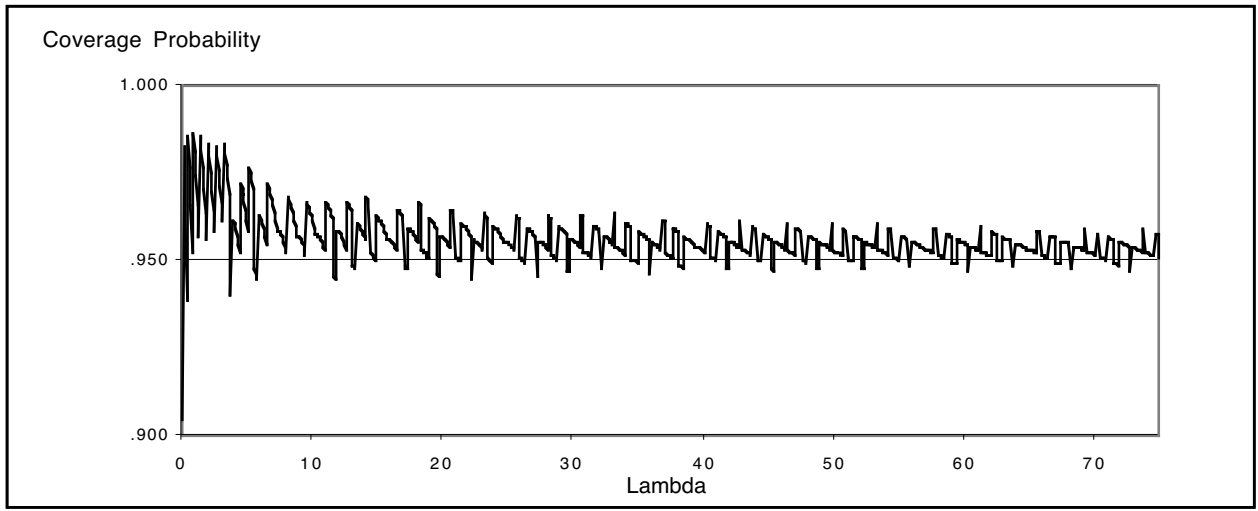

FIGURE 9: Coverage Probabilities for 95% Intervals for $f_1 = 2x + 1$ and $f_2 = 2x + 2$

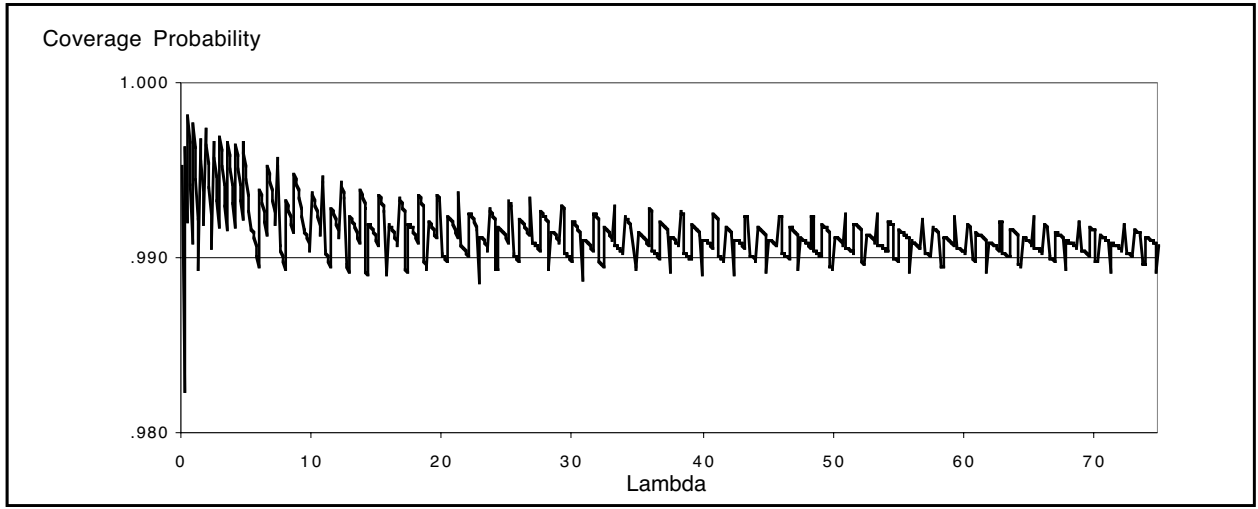

FIGURE 10: Coverage Probabilities for 99% Intervals for $f_1 = 2x + 1$ and $f_2 = 2x + 2$

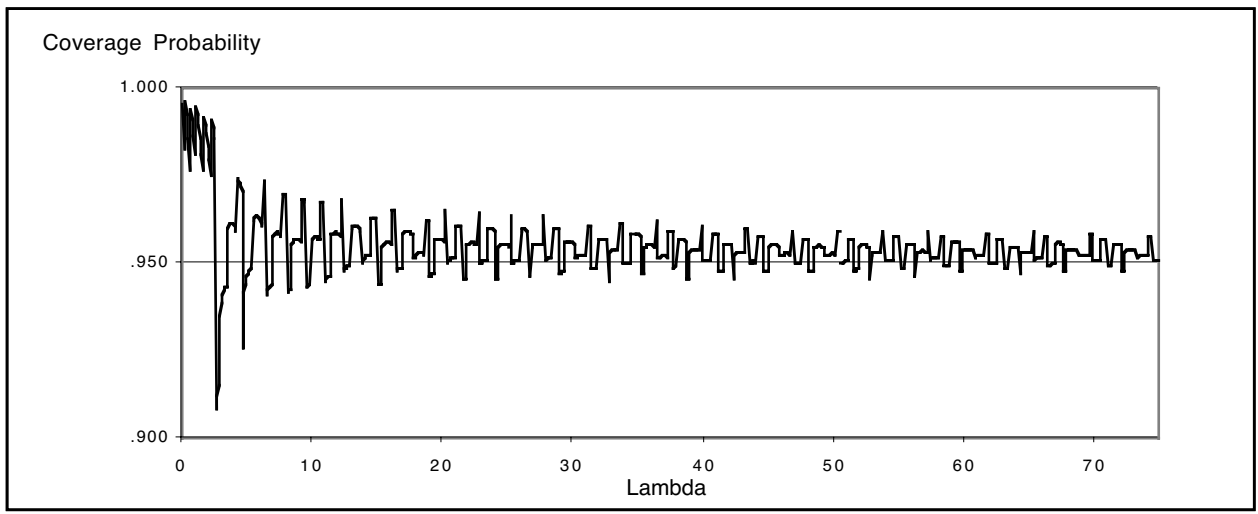

FIGURE 11: Coverage Probabilities for 95% Intervals for $f_1 = 2x$ and $f_2 = 2x + 1$

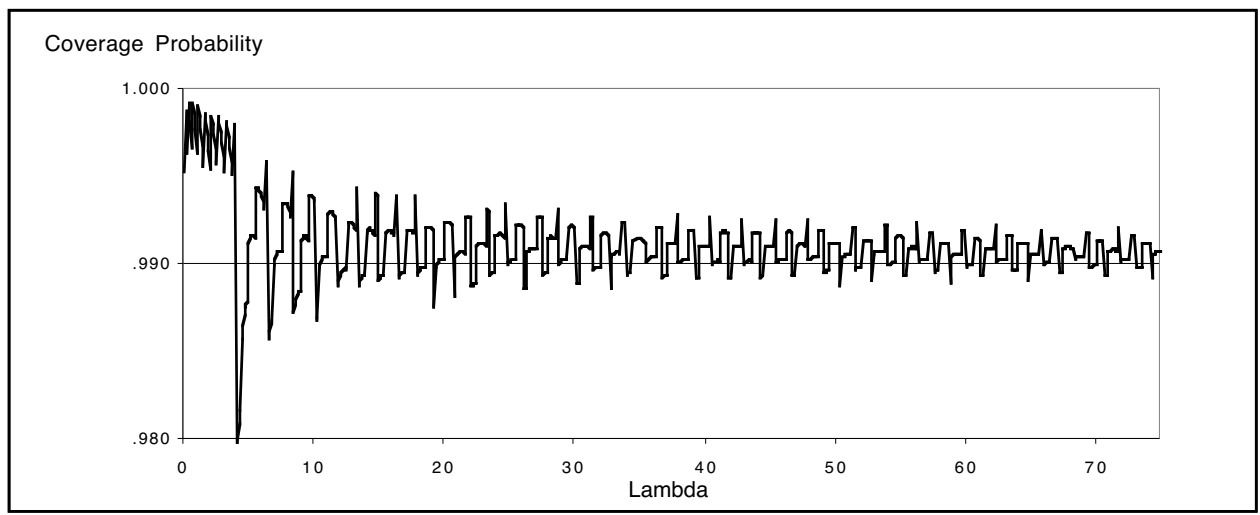

FIGURE 12: Coverage Probabilities for 99% Intervals for $f_1 = 2x$ and $f_2 = 2x + 1$

Figures 9 and 10 show that the case $f_1 = 2x + 1$ and $f_2 = 2x + 2$ produces a very low coverage probability for only a few very small values of $\lambda$.   Except for these few cases, the coverage probability is high for lower values of $\lambda$ and fluctuates around a value somewhat above $1 - \alpha$ for other of $\lambda$ values, crossing occasionally below $1 - \alpha$.

Figures 11 and 12 show that the case $f_1 = 2x$ and $f_2 = 2x + 1$ produces a very low coverage probability for a few small values of $\lambda$.   Except for these few cases, the coverage probability is high for lower values of $\lambda$ and fluctuates around $1 - \alpha$ for other of $\lambda$ values, but is generally somewhat above $1 - \alpha$.   This is similar to the Bayes interval based on the Jeffreys prior for the Poisson.

Either of these intervals might appeal to investigators who prefer approximate estimation intervals rather than "exact" ones to obtain a coverage probability that is close to $(1-\alpha)$ rather than one that is always at least $(1-\alpha)$ but often much larger than $(1-\alpha)$.   On, the other hand, neither of these intervals has a theoretical justification as provided for the classical, structural/Bayes or Bayes intervals explored.

## 6.  Summary statistical measures

The performance of each of the intervals above has been illustrated via a plot of coverage probabilities determined for $\lambda = 0.1$ to $75.0$ in steps of $0.1$.   For a numerical summary, the mean, minimum and maximum coverage probabilities were also determined for the values of $\lambda$ studied.

These summary statistics are tabulated below in Tables 1 and 2.

TABLE 1:Coverage probabilities for 95% Intervals

| Interval Basis | $f_1$ | $f_2$ | Mean | Minimum | Maximum |
| --- | --- | --- | --- | --- | --- |
| "Usual" | $2x$ | $2x+2$ | 0.9611 | 0.9504 | 0.9964 |
| Structural & Bayes-Jeffreys prior for a positive parameter | $2x$ | $2x$ | 0.9473 | 0.8701 | 0.9964 |
| Bayes-uniform prior | $2x+2$ | $2x+2$ | 0.9497 | 0.8187 | 0.9743 |
| Bayes-Jeffreys prior for Poisson | $2x+1$ | $2x+1$ | 0.9499 | 0.9048 | 0.9865 |
| "Other" - option 1 | $2x+1$ | $2x+2$ | 0.9561 | 0.9048 | 0.9865 |
| "Other" - option 2 | $2x$ | $2x+1$ | 0.9549 | 0.9086 | 0.9964 |

Table 2: Coverage probabilities for 99% Intervals

| Interval Basis | $f_1$ | $f_2$ | Mean | Minimum | Maximum |
| --- | --- | --- | --- | --- | --- |
| "Usual" | $2x$ | $2x+2$ | 0.9926 | 0.9902 | 0.9992 |
| Structural & Bayes-Jeffreys prior for a positive parameter | $2x$ | $2x$ | 0.9891 | 0.9653 | 0.9992 |
| Bayes-uniform prior | $2x+2$ | $2x+2$ | 0.9898 | 0.9048 | 0.9947 |
| Bayes-Jeffreys prior for Poisson | $2x+1$ | $2x+1$ | 0.9900 | 0.9736 | 0.9982 |
| "Other" - option 1 | $2x+1$ | $2x+2$ | 0.9915 | 0.9825 | 0.9982 |
| "Other" - option 2 | $2x$ | $2x+1$ | 0.9912 | 0.9788 | 0.9992 |

## 7.  Conclusions

Several chi-square intervals are available as alternatives to the 'standard' confidence interval for a Poisson parameter.  These intervals provide options for experimenters who wish to have shorter intervals than the 'standard' and who may have a preference for a coverage probability that is generally close to the nominal $(1 - \alpha)$, although sometimes less than $(1 - \alpha)$.  Although these intervals are not "exact" (the coverage probability is not always at least $(1 - \alpha)$,) they require no preparation of computer code nor the running of algorithms.  As well, for experimenters who use Bayesian or structural inference methods, chi-square-based estimation intervals with post-sample probability $(1 - \alpha)$ of capturing the true value of the Poisson parameter are available.  Exploring the coverage probability also provides some validation for "objective" priors for Bayesian analysis, particularly the Jeffreys prior particular to the Poisson model.